\documentstyle[12pt]{article}

\newcommand{\beq}{\begin{equation}}
\newcommand{\eeq}{\end{equation}}
\newcommand{\beqcol}{\begin{array}{rcl}}
\newcommand{\eeqcol}{\end{array}}

\newcommand{\g}{{\bf g}}
\newcommand{\h}{{\bf h}}

\begin{document}
{\thispagestyle{empty}
\begin{flushright}  YITP/U-95-32 \\ August 1995 \end{flushright}
\vfill
\begin{center}
{\Large
%{\obeylines
On the deformability of Heisenberg algebras
}
\\
\vspace{2cm}
{\sc Mathias Pillin}               \\
\vspace{1cm}
Uji Research Center             \\
Yukawa Institute for Theoretical Physics  \\
Kyoto University                          \\
Uji, Kyoto 611, JAPAN     \\
\ \                      \\
e-mail: map@yisun1.yukawa.kyoto-u.ac.jp
\medskip

\vfill
{\bf Abstract}
\end{center}
\begin{quote}
Based on the vanishing of the second Hochschild cohomology
group of the enveloping algebra of the Heisenberg algebra
it is shown that differential algebras coming from
quantum groups do not provide a non-trivial deformation of
quantum mechanics. For the case of a $q$-oscillator
there exists a deforming map to the classical algebra.
It is shown that the differential calculus on quantum planes
with involution, i.e. if one works in position-momentum realization,
can be mapped on a $q$-difference calculus on
a commutative real space. Although this calculus leads to an
interesting discretization it is proved that it
can be realized by generators of the undeformed algebra
and does not posess a proper group of global transformations.
\end{quote}

\eject
}

\setcounter{page}{1}

\section{Introduction}

\bigskip
\bigskip

It is known that the deformation of an algebra, either
of Lie or associative type, is connected to its
(Chevalley or Hochschild) cohomology \cite{GERST}.
More precisely, for an algebra
$\g$ the second cohomology group $H^2( \g , \g )$ contains
the information if a non-trivial deformation of it exists or not.
In particular, if $H^2( \g , \g ) = 0 $ then there exits
no non-trivial deformation of $\g$ .

\medskip

This result can readily be applied to the case of quantum
groups \cite{DRI1, JIM}. Here one takes for example a
finite-dimensional semisimple Lie algebra $\g$ and addresses the
question of existence of deformations of its
enveloping algebra ${\cal U}(\g )$. It is well known
that we have non-trivial deformations denoted by ${\cal U}_{h}(\g )$
\footnote{We take for the deformation parameter
$q = e^h >1$ throughout this paper.}
as long one considers ${\cal U}_{h}(\g )$ as being a Hopf algebra
or at least a bialgebra. The non-triviality of this deformation
comes from the fact that
$H^2({\cal U}(\g ),{\cal U}(\g ))_{bialgebra}\simeq \Lambda^2(\g )
\ne 0$ \cite{DRI2}, \cite{KASSEL} ch.18, where $\Lambda(\g )$ denotes
the exterior algebra.

If one in contrast would consider only the algebra part
of ${\cal U}(\g )$ the classical Whitehead lemma applies in
this case. That lemma states that for a finite
dimensional semisimple Lie algebra $\g$ and and finite-dimensional
left-$\g$-module $M$ it holds that:

\beq
 H^1 (\g , M ) = H^2 (\g , M ) = 0
\label{whitehead}
\eeq

This result gives rise to the following:

\medskip

{\bf Theorem 1.1}\cite{DRI2}
{\it There exists a unique morphism
\beq
    \alpha \quad : \quad {\cal U}_h(\g ) \quad
            \rightarrow \quad {\cal U}(\g )[[h]]
\eeq
of topological algebras, such that $ \alpha \equiv id \quad mod \ h$.
}

\bigskip

Practically speaking this result means that there exist
{\it deforming maps} in the sense of \cite{CURT} which connect the
generators of the deformed to the ones of the undeformed algebra
respectively.

\bigskip

It has been shown in the fundamental paper \cite{BAYEN} that
quantum mechanics itself can be understood as a deformation
of classical mechanics. Already in that work the question
of stability of quantum mechanics with respect to
further deformations has been addressed. This question
was answered affirmatively some years ago
\cite{CLOUX} in showing that the second cohomology
group of the Weyl algebra regarded as a bimodule over
itself is zero. This means that (at least in the reasonable
sense of star-products) quantum mechanics cannot be
further deformed.

\bigskip

However, in recent years the possibility of $q$-deformations of
the Heisenberg algebra has been studied extensively within the
context of quantum groups, see e.g. \cite{AWATA,TAR,FEIN,
HAYA,RIEF,ROS,SW} and references therein. The
characteristic relation arising from these studies is of the type:
\beq
     p x - q  x p  = - i
\label{qheis1}
\eeq

In general one has two possiblities of considering these
algebras according to the inequivalent antilinear involutions
one can choose on the algebra (\ref{qheis1}).

The first choice is to take a Bargmann-Fock type conjugation,
i.e. $\bar{x} =-i p$. In
this case (\ref{qheis1}) becomes a $q$-oscillator algebra.
Here the
equivalence with the undeformed case is known, see e.g. \cite{ZACH}.

There exists, however, a second
possiblity which has been addressed in \cite{SW}. In this case
one wants to interpret the generators of (\ref{qheis1}) as
being momentum and position operators respectively.

The aim of this paper is to show that due to the
mentioned rigidity theorem for the Weyl algebra even
the second approach (in any dimensions) does not provide a
true deformation of the Heisenberg algebra.

\bigskip

The plan of the paper is as follows. In section 2 Heisenberg
and Weyl algebras are defined and the Hochschild cohomology
is calculated. The immediate consequences for $q$-oscillator
algebras are explained. We study the algebra (\ref{qheis1})
with position and momentum operators in section 3. It turns
out that this algebra does not have a proper group of global
transformations containing the Weyl group of ordinary
quantum mechanics as it should in order to have a
quantum mechanical interpretation. In section 4 $q$-difference
algebras on almost commutative spaces with involution on
both the coordinates and the $q$-difference operators will
be considered. The existence of this involution is necessary
for having a position-momentum interpretation of the generators.
A uniqueness result for the calculus will
be obtained. The more general case of a $q$-differential
calculus on the real $SO_q(N)$ quantum plane is
addressed in section 5. It will be proved that there
exists a deforming map which provides an
isomorphsim of this $q$-differential calculus to
the $q$-difference calculus on a commutative space.
Finally we summarize and comment on prospects in section 6.

\bigskip
\bigskip

\section{The cohomology of the Weyl algebra}

\bigskip

We denote by  $\h_n$ the Heisenberg algebra which is generated by
the set of generators $p_1,q^1; p_2,q^2;\ldots ; p_n,q^n$ and the
unity together with the relations:
\beq
[q^i, p_j ] = i \delta^i_j, \quad {\rm and} \quad
[q^i, q^j ] = [p_i, p_j ] = 0 .
\label{heis2}
\eeq
The enveloping algebra of $\h_n$ is denoted by $A_n$ and will
be called Weyl algebra. We note that by minor changes in the
definitions
one can also consider the Weyl algebra $A_n(k)$ over an
arbitrary field $k$ of characteristic zero.

It is obvious that elements of the form
$(p_1)^{i_1}(q^1)^{j_1} (p_2)^{i_2}(q^2)^{j_2} \cdots
(p_n)^{i_n}(q^n)^{j_n}$ with $i_1,j_1,$ $\ldots , $ $i_n,j_n \in
{\bf N}_0$ generate $A_n$ as a vector space.

One naturally has a realization of the generators of $A_n$ by
differential operators:
\beq
q^k \sim x^k, \qquad p_k \sim -i {{\partial}\over{\partial x^k}} .
\label{realis}
\eeq

\bigskip

We can now state the following theorem on the cohomology
$ H^{\star} ( A_n (k), A_n (k))$.

\medskip

{\bf Theorem 2.1} \cite{CLOUX}:
{\it It holds that:
\beq
 H^{m} ( A_n (k), A_n (k)) = k \delta_{m,0}
\eeq
}
\bigskip

Due to the importance of this result we sketch the proof.

{\sc Proof}: One first proves a more general result. Let {\bf a} be a
finite-dimensional nilpotent Lie algebra over $k$ and $B$ its
enveloping algebra. We denote by {\bf b} an ideal of
{\bf a} and by $\lambda$ a character of {\bf b}. For
$\xi \in {\bf b}$ the set of elements in $B$ of the form
$\xi - \lambda(\xi )$ is denoted by ${\bf b}_{\lambda}$.
Obviously ${\bf b}_{\lambda}$ is a sub-vector space of
$B$ stable under the adjoint action of {\bf a}. Moreover
we define $ B_{\lambda}:= B / B {\bf b}_{\lambda} $.
Using the {\it inverse process} of homological algebra
(see e.g. \cite{CE}) one can show that for a
$ B_{\lambda}$-bimodule $X$ the following statement is
true:
\beq
H^m ( {\bf a} / {\bf b} , X ) \simeq H^m (  B_{\lambda} ,X) ,
\qquad \forall m  \ge 0 .
\label{cloux1}
\eeq
We can now pass to the special case of the Weyl algebra.
In this case ${\bf a}$ is the Heisenberg algebra ${\bf h}_n$ and
${\bf b}= {\bf z}$ its center which is of course trivial. Obviously
it holds that $B_{\lambda} = A_n(k)$ with $2n+1$ being
the dimension of ${\bf a}$. In order to prove
the theorem we set $X = A_n(k)$. By (\ref{cloux1}) we
are enabled to take $H^{\star} ( A_n(k),A_n(k)) $ $=$
$ H^{\star} ( {\bf h}_n/k {\bf z},  A_n(k) )$.

Using the realization of the generators of ${\bf h}_n$ in terms
of differential operators (\ref{realis}) the problem
of calculating $ H^{\star} ( {\bf h}_n/k {\bf z},  A_n(k) )$
is mapped to a problem in de Rham cohomology.

The observation that all the derivations in $A_n(k)$ are
interior ones completes the proof of the theorem.

$\Box$

\bigskip

The essence of this theorem is that there does not exist a
non-trivial deformation of the Heisenberg algebra within the
category of algebras.

That $q$-deformations of Heisenberg algebras are trivial
from the general deformation theory point of view is evident
if one considers $q$-deformed oscillator algebras. In one
dimension this algebra takes the form:

\beq
a a^{\dagger} - q a^{\dagger} a = 1
\label{qboson}
\eeq

Studying its representations one finds
the same Fock space as in the undeformed case. The only
difference between the classical and the $q$-case are
the norms of the operators.

The undeformed oscillator algebra is given by:

\beq
A A^{\dagger} -  A^{\dagger} A = 1, \qquad N:= A^{\dagger} A .
\label{osc}
\eeq

If one chooses a certain completion of these algebras in
the $h$-adic topology it is possible to introduce an
element of the form $exp{(h N)} = q^N$. In the sense of
theorem 1.1 we then get the following deforming maps:

\beq
a = q^{{N}\over{4}} A \sqrt{ {[N]}\over N} , \quad
a^{\dagger} = q^{{N}\over{4}} \sqrt{ {[N]}\over N } A^{\dagger},
\qquad [N] : = { { q^{N\over 2} - q^{-{N\over 2}}} \over
               {q^{1\over 2}-q^{-{1\over 2}}}} .
\label{osc2}
\eeq

We note that deforming maps are not unique in general.

This deforming map can be applied to higher dimensional
$q$-boson algebras as well. For instance it has been
shown in \cite{OGI1} that the $q$-differential calculus
on quantum planes belonging to quantum groups $SL_q(N)$
and $SO_q(N)$ (see section 5 for more details)
can be mapped into a tensor product of mutually
commuting algebras
of the form (\ref{qboson}). However, this result applies
only to cases without any reality conditions neither
on the quantum plane nor on the differential operators. This
means that the mentioned map in \cite{OGI1} is not at all
compatible with the natural involution (real form) coming
from the quantum group.

\bigskip
\bigskip

\section{Remarks on ``$q$-deformed quantum mechanics''}

\bigskip
\bigskip

As mentioned in the introduction there exists another
approach to the deformation of Heisenberg algebras which
comes from the Wess-Zumino differential calculus \cite{WZ}.
In this approach one considers the quantum plane
which is a certain comodule of a quantum group and the
$q$-differential calculus on it. The differential relations
are interpreted as $q$-Heisenberg relations.

In contrast to the $q$-boson algebras the involution in this
approach is not of Bargmann-Fock type. As mentioned above it
has been shown in \cite{OGI1} that with an involution of Bargmann-Fock type
the full $q$-differential algebra in the case of
$SL_q(N)$ and $SO_q(N)$ can be transformed into a
tensor product of mutually commuting algebras similar to (\ref{qboson}).

\medskip

The authors of \cite{SW} investigate the following algebra:
\beq
    p x - q x p = - i
\label{qheis2}
\eeq

Like in ordinary quantum mechanics the generators $p$ and $x$
ought to be interpreted as momentum and position operators on some
Hilbert space. Since $q$ is taken to be real one cannot
find an antilinear involution which allows for taking
both generators to be real under involution. However, one can
take $\bar{p} = p $ and introduces an additional generator
$\bar{x}$ together with the obvious relations:
\beq
  p \bar{x} - q^{-1} \bar{x} p = - i q^{-1}, \qquad
  x \bar{x} = q  \bar{x} x .
\label{qheis3}
\eeq
For convenience an additional object is introduced by using
usual commutators:
\beq
 r = i [ p , x ] , \qquad  \bar{r} = i [ p , \bar{x} ] .
\label{redef}
\eeq

The aim is to get a $q$-Heisenberg algebra with formal real
objects. If $p$ is interpreted as real momentum an obvious
choice for a real position would be $\tilde{\xi}:= x + \bar{x}$.
This definition results in the algebra:

\beq
  \tilde{\xi} p - q^{-1} p \tilde{\xi} = (q^{-1} +1) i \bar{r} ,
  \qquad
  \tilde{\xi} p - q  p \tilde{\xi} = (q^{-1} +1) i q r .
\label{qheis4}
\eeq

The problem with this algebra is that even though $r$ $q$-commutes
with $p$ it has rather involved relations with $\tilde{\xi}$.
To circumvent this problem one uses the observation that
$r$ and $\bar{r}$ can be decomposed, however non-uniquely,
into a formal real and and a quasi unitary object.

Denoting $T := r \bar{r}$ one is allowed to write
$\bar{r}:= \sqrt{q} u T^{1/2}$ and $ r := \sqrt{q} T^{1/2} \bar{u}$ which
of course implies $u\bar{u} = q^{-1} = \bar{u} u$. Applying
these decompositions to the algebra (\ref{qheis4}) and redefining
the position to be
\beq
    \xi := { \sqrt{q} \over { q +1} }
            \left( T^{-1/2} x + \bar{x} T^{-1/2} \right) ,
\label{xidef}
\eeq

yields the following algebra:
\beq
  \xi p - q^{-1} p \xi =  i u ,
  \quad
  \xi p - q  p \xi =  i u^{-1},
  \qquad
  u p = q^{-1} p u, \quad u \xi = q^{-1} \xi u .
\label{qheis5}
\eeq

Although this algebra suggests to be interpreted as a deformation
of the Heiseberg algebra it turns out that all relations can
be entirely realized within a certain completion of the
enveloping algebra of the usual Heisenberg algebra
$[x_c,p_c] = i$.

One has the freedom to interpret the momentum $p$ to be the usual
momentum generator $p_c$. This gives:

\beq
    u = exp( -i h p_c x_c), \qquad
   \bar{u} = exp( i h x_c p_c) = q^{-1} u^{-1} .
\label{ureal}
\eeq

The following realization follows directly from (\ref{qheis5}):

\beq
    \xi = { i \over p} { {u - u^{-1}}\over{q-q^{-1}}}
        \equiv x_c \quad mod \quad  h
\label{xireal}
\eeq

\bigskip

It is shown in \cite{SW,JOS} that the spectra of $p$ and
$\xi$ are discrete. In momentum representation the
Hilbert space states are given by vectors $ |n \rangle_{\pi_0}$
where $n \in {\bf Z}$ and the continuous parameter
$\pi_0 \in [ 1, q)$ labels the the different irreducible
representations of the algebra (\ref{qheis5}). One then
has :

\beq
 p |n \rangle_{\pi_0}  = \pi_0 q^n |n \rangle_{\pi_0}, \qquad
 \xi |n \rangle_{\pi_0} = { i \over {\pi_0 q^n (q-q^{-1} )}}
             \left( q^{1\over{2}} | n-1 \rangle_{\pi_0} -
            q^{-{1\over{2}}} | n+1 \rangle_{\pi_0} \right)   .
\label{pspec}
\eeq

By a proper Fourier transformation \cite{JOS} it is possible to show
that one can also construct irreducible representations in
which $\xi$ is diagonal. Its spectrum then is similar to
the one of $p$ in (\ref{pspec}). Actually the momentum eigenstates
$ |n \rangle_{\pi_0}$ can be realized by ordinary functions.
Up to normalization we have :
\beq
|n \rangle_{\pi_0}  \sim exp \left( i q^n \pi_0 x_c  \right)
\label{spep}
\eeq

The operator $p=p_c$ then acts on these states by ordinary
Schr\"odinger representation although the application
of $x_c$ does lead out of the irreducible representations
of the algebra (\ref{qheis5}).

Although it seems quite interesting to interpret the
algebra (\ref{qheis5}) to be a $q$-deformation of
the Heisenberg algebra which provides a discretization
quite similar to the one arising from ordinary lattice
quantum mechanics we have the following

\medskip

{\bf Lemma 3.1}:
{\it The algebra (\ref{qheis5}) does not allow for
a group of global transformations which contains the
Weyl group of usual quantum mechanics consistent
with the irreducible representations of that algebra.
}
\medskip

{\sc Proof}: We have seen that all generators
of the algebra (\ref{qheis5}) can be realized in terms of the
generators of the classical Heisenberg algebra.

Strictly
speaking (\ref{qheis5}) is not generated by the classical
Heisenberg algebra itself which is nilpotent but by the
minimal solvable extension of it. This means that on
the algebra level we have an additional generator
equivalent to $p_c x_c$.
The enveloping algebras of these algebras are identical since
the additional generator lies in the vector space spanned
by elements of the form $(p_c)^i ( x_c)^j$ with $i,j \in {\bf N}$.

If we denote the classical nilpotent Heisenberg algebra
by $H_1$ the corresponding group is the Weyl group $W_1$.
Obviously $H_1$ is a subalgebra and moreover a Lie bi-ideal
of its minimal solvable extension.

The usual uniform lattice discretization corresponds to
considering the Weyl group $W_1$ not over the real numbers but
over the integers ${\bf Z}$.

The group corresponding to the minimal solvable extension
of $H_1$ is some product of $W_1$ with a group of dilatations
generated by $p_c x_c$ which we will call $D$. We
denote this fact by $WD_1 = W_1 \bullet D$.

The spectra of (\ref{qheis5}) would correspond -- by
the conformal invariance of the real line -- to fix
the group parameters in the pure $D$ part of $WD_1$
to the set $h{\bf Z}$. For consistency the
space of group parameters for the full group $WD_1$
has to be $h{\bf Z}$. That the global transformations
corresponding to this set is inconsistent with the
spectrum of (\ref{qheis5}) is evident. The only thing
one has to do is to apply a translation of the form
$exp(\alpha p_c)$ to the Hilbert space of the algebra
(\ref{qheis5}) with $\alpha / h \in {\bf Z}$. This
completes our proof.

$\Box$

\bigskip

This result can be applied to higher dimensional cases.
In \cite{HEB} for example the case of $SO_q(N)$ covariant
quantum mechanics has been investigated. It turned out
that the total Hilbert space of the theory is a tensor product
of two Hilbert spaces. One of them corresponds to ${\cal U}_q(su2)$
and the other is the one coming from (\ref{qheis5}).
In the sense of theorem 1.1 and of lemma 3.1 that models
do not provide a true deformation of quantum mechanics.

\bigskip
\bigskip

\section{q-derivatives on almost commutative spaces}

\bigskip

In this section $q$-difference calculi on commutative or
almost commutative spaces with involution are considered.
By almost commutative (sometimes also referred to as $q$-commutative)
we mean two objects $A$ and $B$ having
a commutation relation of the form $AB = q^r B A$ where $r$
might be any number different from zero. Having the background
of quantum groups we require that for $q=1$ we obtain the
usual continuous calculus.

\medskip

Let us first assume that the configuration space is a
finite-dimensional commutative
algebra generated by objects $x^{\alpha}$. We will talk about
these generators freely as coordinates. For the basis we
choose the following convention. If the dimension of the space
is odd (dim=$2n+1$) then
$\alpha \in {\cal I}_{2n+1} :=
\{ -n,\ldots, -1, 0, 1, \ldots , n \}$ if it is even (dim =$2n$)
we have that  $\alpha \in  {\cal I}_{2n} :=
\{ -n, \ldots, -1, 1, \ldots , n \} $. When writing $\alpha \in
{\cal{I}}$ it is assumed that $\alpha$ is either in ${\cal I}_{2n+1}$
or in ${\cal I}_{2n}$.

An (anti-linear) involution on this space is introduced by the
rule $ \overline{x^{\alpha}} := x^{-\alpha}$. This
is just the standard involution on an euclidean space in a
lightcone like basis. We remark that our results are
of course basis-independent. Since we later want to make
contact to the quantum group case the chosen basis is
convenient.

Let $D_{\alpha}$ be $q$-partial derivatives acting on $x^{\alpha}$.
For the application we have in mind and hence
in analogy to (\ref{qheis5}) we require these
partial derivatives to have the standard conjugation property
$\overline{D_{\alpha}} = - D_{ - \alpha}$.

As has been mentioned in the introduction the existence of an involution
of this kind is essential for our considerations because
in the sense of the previous section we have in mind a
position-momentum interpretation of the differential algebra.

\smallskip

For the remainder of this paper no summation over repeated indices
is assumed.

\smallskip

By ${\rm Diff}_{q^{k(\alpha )}}(x^{\alpha}) $ we denote the algebra
generated by $x^{\alpha}, D_{\alpha},  u_{\alpha} $ and
$u_{\alpha}^{-1}$ where $k(\alpha )$ is some number not
equal to zero belonging
to the index $\alpha$. The ideal of relations in
${\rm Diff}_{q^{k(\alpha )}}(x^{\alpha}) $ is generated by:

\beq
\beqcol
D_{\alpha} x^{\alpha} -  q^{k(\alpha)} x^{\alpha} D_{\alpha}& = &
u_{\alpha}^{-k(\alpha)},  \\
D_{\alpha} x^{\alpha} -  q^{-k(\alpha)} x^{\alpha} D_{\alpha}& = &
u_{\alpha}^{k(\alpha)},
\eeqcol \qquad
\beqcol
u_{\alpha} x^{\alpha} &=& q x^{\alpha} u_{\alpha} ,\\
u_{\alpha} D_{\alpha} &=& q^{-1} D_{\alpha} u_{\alpha}
\eeqcol
\label{diffdef}
\eeq

\bigskip

{\bf Proposition 4.1}
{\it Under the previous assumptions together
with $k(\alpha) = -k(-\alpha) $ and $\alpha \in {\cal{I}}$
it holds that:

Every linear $q$-difference calculus on a commutative space,
obeying the formal hermiticity conditions $\overline{x^{\alpha}}
= x^{-\alpha}$ and $ \overline{ iD_{\alpha}} = iD_{-\alpha}$
is equivalent to
\beq
  \bigotimes_{\alpha \in \cal{I}} {\rm Diff}_{q^{k(\alpha )}}(x^{\alpha}).
\label{propform}
\eeq
}
\bigskip

The tensor product here implies that all the off-diagonal relations
are commutative.

Before proving the proposition we state the following

\medskip

{\bf Corollary 4.2}
{\it The usual continuous partial derivatives on the coordinates
are given by $\partial_{\alpha} x^{\beta} - x^{\beta}
\partial_{\alpha} = \delta^{\beta}_{\alpha}$. The formal hermiticity
of the relations (\ref{diffdef}) restricts $u_{\alpha}$ to be of the
form $ u_{\alpha} = exp\left( h \ x^{\alpha} \partial_{\alpha}
\right)$. Moreover we have $\overline{u_{\alpha}} =
q^{-1} u_{-\alpha}^{-1}$.
}
\medskip

{\sc Proof}: If we assume that the proposition is true, i.e.
relations (\ref{diffdef}) hold for all $\alpha \in {\cal{I}}$,
the corollary is a direct consequence of these relations.

Since we require that for $q \rightarrow 1 $ the
derivatives $D$ should turn into the continuous derivatives it
must hold for all $\alpha$ maybe up to some normalization that
$D \equiv \partial$ $mod$ $h$.

Together with the linearity of the calculus and the formal
hermiticity of the $q$-derivatives this requires
that in the realization of $D_{\alpha}$ some linear combination of
the $ u_{\alpha}^{k(\alpha)}$ and $ ( u_{\alpha}^{k(\alpha)} )^{-1}$
must appear.

Hence, a realization of any $D_{\alpha}$ is of the form:
\beq
 D_{\alpha} \sim {1 \over x^{\alpha} } {{ u_{\alpha}^{k(\alpha)} -
   \left(u_{\alpha}^{k(\alpha)} \right)^{-1}} \over {
    q^{ k(\alpha)} - q^{- k(\alpha)}}}
   \prod_{\beta} u_{\beta}^{r(\beta)} + \lambda I
\eeq

The product appearing in this expression is taken over
some $u$'s such that the hermiticity of $iD$ is not
spoiled. $\lambda$ is a number tending to zero as
$ q \rightarrow 1$. $I$ is a polynomial in
$x$'s,$D$'s, and maybe classical $\partial$'s and is required to have
the dimension of a partial derivative.

It is easy to show that the $q$-derivative in the above
expression can be shifted and rescaled. This means
that according to the assumptions of the proposition
$D_{\alpha} $ has the following realization:

\beq
D_{\alpha}  = {1 \over x^{\alpha} } {{ u_{\alpha}^{k(\alpha)} -
   \left(u_{\alpha}^{k(\alpha)} \right)^{-1}} \over {
    q^{ k(\alpha)} - q^{- k(\alpha)}}}
\eeq

Now it is easy to conclude. The only thing which has to
be done is to calculate the commutation relations of the
so realized $D_{\alpha}$ with all other generators. The
relation with $x^{\alpha}$ is just (\ref{diffdef}) while
it commutes with $D_{\beta}$ and $x^{\beta}$ for
$\alpha \ne \beta$.

$\Box$

\bigskip

We note that the considerations of section 3 apply to each
component of the tensor product (\ref{propform}).

\medskip

{\it Remark 1}: It can easily be seen that if the
coordinate algebra is almost commutative the corresponding
$q$-difference calculus can be transformed into the
form described in proposition 4.1.

To illustrate this we use a simple three-dimensional
example which, however, is generic. We take
commutative coordinates $x^{-1}$, $
x^{0}$ and $x^{1}$ subject to the above mentioned
conjugation rule. We can make these coordintes
$q$-commuting by using the $q$-shift operators $u_{\alpha}$.
In the simplest case we set:

\beq
 x^{-1}_q := x^{-1}, \quad    x^{0}_q := u_{-1}(u_{1})^{-1} x^{0},
  \quad    x^{1}_q := x^{1} .
\label{qcoord}
\eeq
By corollary 4.2 this unitary transformation preserves the
real structure of the coordinates and leads to the
commutation relations:

\beq
 x^{-1}_q x^{0}_q = q  x^{0}_q  x^{-1}_q , \quad
 x^{0}_q x^{-1}_q = q^{-1}  x^{0}_q   x^{1}_q  \quad
 x^{-1}_q x^{1}_q = x^{1}_qx^{-1}_q  .
\label{commutrel}
\eeq

The $q$-derivative with respect to $x^{0}$ has to be rescaled
as well by $D_0^q := (u_{-1}(u_{1})^{-1})^{-1} D_0$ while the
other ones remain the same. The outcome is a $q$-difference
calculus in which the diagonal relations are almost the same
as (\ref{diffdef}) while the off-diagonal ones are
$q$-commutative.

Although the Hilbert spaces of the commutative and of the
almost commutative case are not identical the algebras
can simply be related.

These almost commutative calculi appear for example by reduction
of a $GL_q$-quantum group to some lower dimensional
orthogonal quantum group (see e.g. \cite{CAST}).

\bigskip

{\it Remark 2}: It is a well known problem in the
study of inhomogeneous quantum groups (e.g. \cite{OGI3}) that it
is difficult to find a coproduct which preserves
the formal antihermiticity of $q$-derivatives. It
can be read off the algebra (\ref{diffdef}) that
the comultiplication of the formal antihermitian
$q$-difference operator necessarily involves either
the quantity $u$ or $u^{-1}$. Due to the
conjugation property stated in corollary 4.2 it
is clear that the comultiplication of a
$q$-difference operator can hardly preseve
the formal anti-hermiticity.

\bigskip

{\it Remark 3}: Another problem with the differential
calculus on quantum spaces is the nonlinear conjugation
rule of the derivatives \cite{OGI2}. Although we will
treat this case in the next section in some detail,
already at this stage some comments are necessary.

The $q$-partial derivatives as they turn out of the
Wess-Zumino calculus \cite{WZ} are unsymmetric in the
sense that they produce $q$-shifts in only one
direction in contrast to the case considered in
(\ref{pspec}). It can be shown that
the nonlinear conjugation rule occurs already in
simpler cases.

Let us use the commutative three dimensional space
introduced in remark 1, and introduce unsymmetric
but commuting $q$-difference operators by:

\beq
\partial_{-1} := { u_0 \over{ x^{-1}}}
                 { { 1 - (u_{-1})^2} \over {1-q^2}} , \quad
\partial_{0}  := {{u_{-1} u_{1} }\over{x^{0}}}
                 { { 1 - (u_{0})} \over {1-q}} ,       \quad
\partial_{1}  :=  {u_{0}\over{x^{1}}}
                 { { 1 - (u_{1})^2} \over {1-q^2}} .
\label{part}
\eeq

The action of these derivatives on the coordinates looks
almost like the ones coming from the $SO_q(3)$ covariant
calculus. The diagonal relations are for instance of the form
$  \partial_{0} x^{0} = u_{-1} u_{1} + q x^{0} \partial_{0}$.
The off diagonal relations are almost commutative. If one
introduces a formal antilinear involution, denoted as
above by a bar, on the derivatives
we get e.g.
$ \hat{\partial}_{0} x^{0} = (u_{-1} u_{1})^{-1}
 + q^{-1} x^{0} \hat{\partial}_{0}$ where
$\hat{\partial}_{0}:= -q^{-3} \overline{\partial_{0}}$. The
relation between $\hat{\partial}_{0}$ and ${\partial}_{0}$
is established introducing the quantity
$\Lambda:= u_{-1}^2  u_{0}^2 u_{1}^2 $. We get:

\beq
\hat{\partial}_{0}  \sim \Lambda^{-1} \partial_0
            -{{q+1}\over{q^{-2}-1}} { {(u_{-1} u_{1})^{-1}}\over{x^0}}
	    \left( (u_0)^{-2} -1 \right)
\label{conjr}
\eeq

Analogous relations hold also for the other $q$-partial
derivatives. Thus, the aim of this remark is that the
nonlinear conjugation rule is not surprising when
considering unsymmetric $q$-difference operators in the sense
of (\ref{part}).

\bigskip
\bigskip

\section{q-differential algebras coming from orthogonal quantum
groups}

\bigskip

In this section we want to extend the results obtained in the
previous section to the case of differential calculi on
quantum planes coming from orthogonal quantum groups. These
quantum planes fit well in our treatment since the
compact form of an orthogonal quantum group $SO_q(N)$
naturally induces a real structure on the corresponding
quantum plane.

In order to present the main result of this section we
define the ingredients and fix the notation. For further
details see \cite{OGI2} and references therein.

\medskip

The $\hat{R}$-matrix for the orthogonal quantum group
in $N$ dimensions, $SO_q(N)$, posesses a decomposition
into the following projection operators:
the $q$-analogs of the symmetrizer $P^{+}$, the antisymmetrizer
$P^{-}$, and the trace projector $P^{0}$. The latter defines
the $q$-analogue of the metric tensor $g_{ij}$ by
$P^{0\ ij}_{\ \ \ \ kl} = c g^{ij}g_{kl}$, where $c$ is
some constant. We then have:

\beq
\hat{R} = q P^{+} - q^{-1} P^{-} + q^{1-N} P^{0}
\label{Rdec}
\eeq

The quantum plane corresponding to $SO_q(N)$ is the
algebra generated by $N$ generators $x^{i}$ with $i$
running through the index sets ${\cal{I}}_{2N+1}$ or
${\cal{I}}_{2N}$ which have been defined at the beginning
of the previous section.

The ideal of relations in this algebra is generated by
$\sum_{k,l} P^{-\ ij}_{\ \ \ \ kl}x^k x^l = 0$. The so defined
algebra will be denoted by $V_q(N)$. By definition
$V_q(N)$ is a $SO_q(N)$-comodule. The metric defines
a $SO_q(N)$-invariant object (by the comodule mapping)
$L:= (1+q^{N-2})^{-1}( \sum_{i,j}g_{ij} x^i x^j )$
which is central in the algebra $V_q(N)$.

Due to the real form of the quantum group there exists an
antilinear involution on $V_q(N)$ given by:

\beq
  \overline{x^i} = \sum_j g_{ji} x^j
\label{reali}
\eeq

It has been shown in \cite{WZ} that it is possible to
construct an algebra of $q$-partial derivatives on $V_q(N)$.
We denote these derivatives by $\partial_i$ with
$i \in {\cal{I}}$.
The set of $\partial_i$ spans a $SO_q(N)$-comodule algebra
(just another quantum plane)
with relations dual to the ones in $V_q(N)$, namely:
\beq
 \sum_{i,j} P^{-\ ij}_{\ \ \ \ kl}\partial_i \partial_j = 0
\label{ablrel}
\eeq
The element $\Delta := (1+q^{N-2})^{-1}( \sum_{i,j}g^{ij} \partial_i
\partial_j )$ is central in the algebra of the $q$-derivatives
and invariant under $SO_q(N)$-coaction.

The action of the $q$-partial derivatives on the
generators of $V_q(N)$ is given by:

\beq
\partial_i x^j = \delta^j_i + q \sum_{k,l} \hat{R}^{jk}_{il}x^l \partial_k
\label{action1}
\eeq

This action is unsymmetric by its definition in the sense
of remark 3 in the previous section.

\medskip

The involution on $ V_q(N)$ as defined in (\ref{reali}) cannot
be extended to the $q$-differential algebra. Another copy
of partial derivatives has to be introduced. Using the
notation of (\ref{reali}) we have
$\hat{\partial}_i := -q^N \sum_{k,l}g_{ik} g^{kl} \overline{\partial^l}$.
Although the algebra of these derivatives is generated by
relations identical to (\ref{ablrel}) the action of the
conjugated derivatives is given by:

\beq
\hat{\partial}_i x^j = \delta^j_i + q^{-1}
       \sum_{k,l}\hat{R}^{-1 \ jk}_{\ \ \ \ il}x^l \hat{\partial}_k
\label{action2}
\eeq

For consistency the $q$-derivatives have to satisfy the
relation $ \partial \hat{\partial} = q^{-1} \hat{R}^{-1}
\hat{\partial} \partial$.

The both kinds of $q$-derivatives can be connected using an
element which is similar to the one introduced in remark 3 of
the previous section. Using the classical analogue of
$V_q(N)$ with commutative generators $x^i_{c}$ and the usual
continuous partial derivatives $\partial^{c}_i$
($i \in {\cal{I}}$) we define a classical Euler element
$ E_{c} := \sum_{i \in {\cal{I}}} x^i_{c}\partial^{c}_i$.
This definition gives rise to introduce:

\beq
\Lambda := exp(2h E_{c}), \qquad
        \overline{ \Lambda } = q^{-2N} \Lambda^{-1} .
\label{lamdef}
\eeq

This element can be expressed in terms of the generators
of $V_q(N)$ and the partial derivatives \cite{OGI2}.
A straightforward calculation gives the following almost
commutative relations:

\beq
\Lambda x^k_{(c)} = q^2 x^k_{(c)} \Lambda , \qquad
\Lambda \partial_{k}^{(c)} = q^{-2} \partial_{k}^{(c)} \Lambda.
\label{lamact}
\eeq

The notation $(c)$ means that these commutation relations
hold for both the generators coming from the quantum groups
and for their classical analogs. (\ref{lamact}) shows that
the element $\Lambda$ is well defined in the $q$-algebra and
in the classical algebra as well. Its coproduct is grouplike
but not consistent with the involution on the quantum group
\cite{OGI3}.

Using this element the original and the conjugated derivatives
can be related by the formula:

\beq
 \hat{\partial}_k = \Lambda^{-1} \left( \partial_k +
        q^{N-1} ( q- q^{-1} ) x_k \Delta \right)
\label{qconjr}
\eeq

This equation should be compared with (\ref{conjr}) in the previous
section.

We denote the full algebra generated by $\{ x^i, \partial_i,
\hat{\partial}_i, \Lambda, \Lambda^{-1} | i \in {\cal{I}} \}$
together with their algebraic relations by ${\rm Diff}_{SO_q}^{c}(N)$.

\medskip

In the spirit of section 3 and 4 the task is to construct a
differential calculus on $V_q(N)$ consisting of formal
anti-hermitian $q$-derivatives. The immediate guess for
an object posessing this property is (cf. \cite{FIORE}):

\beq
D_i = \partial_i + q^{-N} \hat{\partial}_i, \qquad
\overline{D_i} = - D_{-i} .
\label{qconjr2}
\eeq

The relations of these newly introduced objects $D_i$ are identical
to (\ref{ablrel}). $D_i$ is a linear combination of
$q$-derivatives acting either via $\hat{R}$ (\ref{action1}) or
via $\hat{R}^{-1}$ (\ref{action2}) on the generators of $V_q(N)$.
Therefore one has two possibilities for writing the diagonal
actions.

\beq
D_i x^i - q^2 x^i D_i = r_i , \qquad
    D_i x^i - q^{-2} x^i D_i = \tilde{r}_i .
\label{diffrelso}
\eeq
This relation is valid for any $i \in {\cal{I}}$ except
for $x^0$ in the odd dimensional case due to the
properties of the $\hat{R}$-matrix. $x^0$ requires
a $q$ rather than a $q^2$ in the diagonal relation. This
should be noted for all the remaining formulae in this
section. The objects $r_i$ and $\tilde{r}_i$ can be calculated
without any difficulties using (\ref{action1}) and
(\ref{action2}) respectively.

The off-diagonal commutation relations are involved.
Moreover the algebra of $r_i$ and $\tilde{r}_i$ with
all generators of ${\rm Diff}_{SO_q}^{c}(N)$  is quite
complicated but their explicit
form is not needed for our purposes.

Due to the properties of the involution (\ref{reali}) and
(\ref{qconjr}) one finds:

\beq
q^{-2} \overline{r_{-i}} = \tilde{r}_i, \qquad \forall i \in {\cal{I}}.
\label{r1}
\eeq

We now proceed in analogy to the one-dimensional case
in section 3 and split the objects $r_i$. For the
decomposition we have in mind the quantities $u_i$
as they have been obtained in corollary 4.2 are
required. Of course, the so defined $u_i$ is not
supposed to obey any managable commutation relation
with the generators of ${\rm Diff}_{SO_q}^{c}(N)$.

Introducing for any $i \in {\cal{I}}$ elements $\rho_i$
we make the ansatz:
\beq
r_i := u_i^{-2} \rho_i  \quad  \Rightarrow \quad
 \tilde{r}_i = \overline{\rho_{-i}} u_i^2
\label{ansatz}
\eeq
The second expression in (\ref{ansatz}) is a consequence
of (\ref{r1}). It is clear that this decompostion is
not unique. However, its existence is guaranteed. Note that
in the odd dimensional case it holds for the zero component
$r_0 = u_0^{-1} \rho_0$.

This ansatz can now be inserted in (\ref{diffrelso})
yielding:
\beq
u_i^{-2} = \left( D_i x^i - q^2 x^i D_i \right)  \rho_i^{-1} , \qquad
u_i^{2}  = \overline{\rho_{-i}}^{-1}
      \left( D_i x^i - q^{-2} x^i D_i  \right) .
\label{udefneu}
\eeq

Now we can apply proposition 4.1 and equation (\ref{diffdef}).
The statement there is that the quantities $u_i$, $u^{-1}_i$,
or any power of them are defined by by $q$-difference relations
on commutative spaces. Thus the factors $\rho_i^{-1}$ and
$\overline{\rho_{-i}}^{-1}$ in (\ref{udefneu}) provide deforming maps
and we have the following

\bigskip

{\bf Theorem 5.1}
{\it The $q$-differential algebra ${\rm Diff}_{SO_q}^{c}(N)$ is isomorphic
to the $q$-difference calculus of the same dimension
$\bigotimes_{\alpha \in \cal{I}} {\rm Diff}_{q^{k(i )}}(x^{i})$
with $k(i) = 2$ for $i > 0$, $k(i) = -2$  for $i<0$, and $k(0) = 1$
as outlined in proposition 4.1.
}

$\Box$

\medskip

We note that this result is not covered by the treatment
in \cite{OGI1}. As mentioned above in the differential
calculus on a commutative space obtained in that paper
is not consistent with an involution on the coordinate
algebra. In contrast that calculus is only consistent
with a Bargmann-Fock type conjugation rule on the
full differential algebra.

\bigskip
\bigskip

\section{Conclusions}

\bigskip

We have pointed out that coming from the differential
calculus on quantum groups one is led to two a priori
non-equivalent approaches towards a $q$-deformation
of quantum mechanics. The first one is based on
$q$-oscillator algebras. The second approach
arises from the $q$-differential calculus on
involutive quantum planes and focusses on
position-momentum interpretation of the generators
of the $q$-differential algebra.

Although both approaches might have interesting
applications, e.g. either in statistical mechanics
or in the theory of generalized hypergeometric series,
it has been shown that due to the rigidity theorem 2.1
for the Heisenberg algebra both approaches do not
yield a true deformation of quantum mechanics.

Moreover it has been shown that it is hardly possible
to find a hermiticity perserving comultiplication
on the generators of both the $q$-difference
calculus and of the $q$-differential calculus.

For these reasons one has to think carefully if
quantum planes, although they have interesting features
\cite{MAP1,MAP2}, provide a reasonable base space
for quantum field theories.

\bigskip
\bigskip

{\bf Acknowledgement}

\bigskip

It is a pleasure to thank R. Sasaki and M. Flato for their
interest in this work and for many intersting and
stimulating discussions. Useful conversations with
D. Sternheimer, H. Awata, and J. Ding are as well acknowledged.

This work was supported by a fellowship of the Japan Society for the
Promotion of Science (JSPS) and by the
Alexander von Humboldt-Gesellschaft.

\end{document}